\documentclass [a4paper,english]{article}
\usepackage{amsmath, amsfonts, amssymb}
\usepackage[cp1251]{inputenc}
\usepackage[T2A]{fontenc}
\usepackage{babel,indentfirst}
\usepackage{graphicx,graphics}

\newcommand{\ba}[2]{\begin{eqnarray}\label{#1} {#2}\end{eqarray}}

\renewcommand{\ge}{\geqslant}

\title{On the correctness of the  main theorem for absolutely  monotonic functions
}
\author{Sitnik S.M.}
\date{}
\begin{document}
\maketitle
\begin{abstract}
\begin{center}
We prove that  the main theorem for absolutely  monotonic functions on $(0,\infty)$
from the book Mitrinovi\'{c}~D.S., Pe\v{c}ari\'{c}~J.E., Fink~A.M. "Classical And New Inequalities
 In Analysis", Kluwer Academic Publishers, 1993, is not valid without additional restrictions.
\end{center}
\end{abstract}

In the classical book \cite{MPF}, chapter XIII, page 365, there is a definition of absolutely  monotonic on $(0,\infty)$ functions.

\textbf{Definition.} A function $f(x)$ is said to be \textit{absolutely  monotonic on $(0,\infty)$} if it has derivatives of all orders and
\begin{equation}
\label{f1}
f^{(k)}(x)\ge 0, x\in (0,\infty), k=0,1,2,\dots \ .
\end{equation}

For absolutely  monotonic functions the next integral representation is essential:

\begin{equation}
\label{f2}
f(x)=\int_0^\infty  e^{xt}\ d\sigma(t),
\end{equation}
where $\sigma(t)$ is bounded and nondecreasing and the integral converges for all $x\in (0,\infty)$.

Also the basic set of inequalities is considered.

 Let $f(x)$ be an absolutely  monotonic function on $(0,\infty)$.  Then
\begin{equation}
\label{f3}
f^{(k)}(x)f^{(k+2)}(x)\ge \left(f^{(k+1)}(x) \right)^2, k=0,1,2,\dots \ .
\end{equation}

After this definition in the book  \cite{MPF} the result which we classify as the main theorem for absolutely  monotonic functions on $(0,\infty)$ is formulated (theorem 1, page 366).

\textbf{The main theorem for absolutely  monotonic functions.}

\textit{The above definition (\ref{f1}), integral representation (\ref{f2}) and basic set of inequalities (\ref{f3}) are equivalent.}

It means:

\begin{equation}
\label{eq}
(\ref{f1}) \Leftrightarrow (\ref{f2}) \Leftrightarrow (\ref{f3}).
\end{equation}

In the book  \cite{MPF} for $(\ref{f1}) \Leftrightarrow (\ref{f2})$ the reference is given to \cite{Wid1}, and an equivalence  $(\ref{f2}) \Leftrightarrow (\ref{f3})$ is proved, it is in fact a consequence of Chebyschev inequality.

In this note we consider a counterexample to the equivalence $(1) \Leftrightarrow (2)$ of Widder. So unfortunately it seems that the main theorem for absolutely  monotonic functions
in the book  \cite{MPF} \textbf{is not valid !!!}.

This counterexample is very simple so it is strange enough it was not found before (cf. also \cite{S}).

Really, consider a function $f(x)=x^2+1$. Obviously for all $x\in [0,\infty)$
\begin{equation}
f(x)\ge 0, f'(x)=2x \ge 0, f''(x)=2 \ge 0, f^{(k)}(x)=0 \ge 0, k>2.
\end{equation}
So this function $f(x)$ is in the class of absolutely  monotonic functions on $(0,\infty)$ due to the definition (1).
If $(1)\Rightarrow (3)$ is valid then the next inequality must be true as a special case of (3) for all $x\in (0,\infty)$

\begin{equation*}
f(x)f''(x) \ge{ \left( f'(x)\right)}^2 \Leftrightarrow 2(x^2+1) \ge 4x^2  \Leftrightarrow 1\ge x^2
\end{equation*}
but this is not valid for all $x\in (0,\infty)$.

As a conclusion we see that implication  $(1) \Rightarrow (3)$ in \cite{MPF} is not valid. It also means that implication $(1) \Rightarrow (2)$ is also not valid. The implication $(2) \Rightarrow (3)$ is obviously valid due to the Chebyschev inequality.

And consequently also the theorem 2 in \cite{MPF}, pages 366--367 on determinant   inequalities  is not valid too if based only on definition (1).

In some papers the above implications are used to derive new results for absolutely  monotonic functions. It seems not to be  a correct way of reasoning.  One way is to change the main theorem   on absolutely  monotonic functions to a proper one, otherwise for all special cases an integral representation must be proved independently.

Comment 1. On the other hand  everything is OK with theorems on completely monotonic functions. An integral representation for them in \cite{Wid2} include the additional condition
\begin{equation*}
\lim_{x\to \infty} f(x)=0.
\end{equation*}
This condition is omitted in \cite{Wid1} but mysteriously mentioned in  \cite{Wid2} with the reference again to
\cite{Wid1}. May be something like it is needed also for absolutely  monotonic functions.

Different aspects of completely monotonic functions are considered in (\cite{Wid1}--\cite{Wid2}), and also for example in the classical expository articles (\cite{Sam1}--\cite{Berg2}).

Comment 2. There are many ways to generalize notions of  absolutely and  completely monotonic functions. It seems that a first step was done by Sergei Bernstein \cite{B1} and very important  generalizations were investigated by Bulgarian mathematicians Nikola Obreshkov \cite{Obr1}--\cite{Obr3} (also known for two celebrated named formulas: Obreshkov generalized Taylor expansion formula and the Obreshkov integral transform) and Jaroslav Tagamlitskii \cite{Tag1}.

Comment 3. With absolute and  complete monotonicity different functional classes are deeply connected: Stieltjes, Pick, Bernstein, Schoenberg, Schur and others, cf. \cite{Berg3}. Just mention  recent papers \cite{Lam1}--\cite{Karp4}.

So the next problems seem to be rather interesting and important.

\textbf{Problem 1.} Give a correct proof for the theorem under consideration from \cite{MPF} and so give justification for equivalences (\ref{eq}).

\textbf{Problem 2.} Generalize the theorem under consideration from \cite{MPF} \textsl{for fractional derivatives} and  give justification for equivalences (\ref{eq}) for this case.

There are applications of considered inequalities in the theory of transmutation operators for estimating transmutation kernels and norms (\cite{Sit3}--\cite{Sit5}) and for problems of function expansions by systems of integer shifts of Gaussians  (\cite{Sit6}--\cite{Sit8}).

The author is thankful to Prof. Ivan Dimovski for useful information on results of N.~Obreshkov and J.~Tagamlitskii.

\end{document}